\documentclass{amsart}
\usepackage{amssymb,amsfonts,amsmath,amsthm,amscd,stmaryrd,dsfont,esint,hyperref,tikz,upgreek,xcolor}
\usepackage[title]{appendix}
\usepackage[mathcal]{euscript}

\usepackage{fullpage}

\theoremstyle{plain}
\newtheorem{thm}{Theorem}
\newtheorem{cor}{Corollary}
\newtheorem{lem}[cor]{Lemma}
\newtheorem{prop}[cor]{Proposition}

\theoremstyle{definition}
\newtheorem{defn}[cor]{Definition}

\usepackage[style=numeric-comp]{biblatex}
\let\div\undefined{}
\let\P\undefined{}
\DeclareMathOperator{\div}{div}

\DeclareMathOperator{\dist}{dist}
\DeclareMathOperator{\cl}{cl}
\DeclareMathOperator{\intr}{int}
\DeclareMathOperator{\disc}{disc}
\DeclareMathOperator{\err}{err}
\DeclareMathOperator{\len}{length}

\newcommand{\weakstarto}{\overset{\ast}{\rightharpoonup}}
\newcommand{\R}{\mathbb{R}}
\newcommand{\Z}{\mathbb{Z}}
\newcommand{\N}{\mathbb{N}}
\newcommand{\Q}{\mathbb{Q}}
\newcommand{\E}{\mathbb{E}}
\newcommand{\P}{\mathbb{P}}
\renewcommand{\tilde}{\widetilde}

\numberwithin{equation}{section}

\title{Quantitative stochastic homogenization of the G equation}
\author{William Cooperman}
\begin{document}
\begin{abstract}
    We prove a quantitative rate of homogenization for the G equation in a random environment with finite range of dependence. Using ideas from percolation theory, the proof bootstraps a result of Cardaliaguet\textendash{}Souganidis, who proved qualitative homogenization in a more general ergodic environment.
\end{abstract}
\maketitle
\section{Introduction}
We consider the behavior, as $\varepsilon \to 0^+$, of the family ${\{u_\varepsilon\}}_{\varepsilon > 0}$ of solutions to the G equation,
\begin{equation}\label{eps-g-eqn}
    \begin{cases}
        D_t u^\varepsilon(t, x) - |D_x u^\varepsilon(t, x)| + V(\varepsilon^{-1}x) \cdot D_x u^\varepsilon(t, x) = 0 &\quad \text{in $\R_{> 0} \times \R^d$}\\
        u^\varepsilon(0, x) = u_0(x) &\quad \text{in $\R^d$},
    \end{cases}
\end{equation}
where $d \geq 2$, $V \in C^{1,1}(\R^d; \R^d)$ is a divergence-free vector field and $u_0 \colon \R^d \to \R^d$ is Lipschitz. The level sets of $u^\varepsilon$ model a flame front which expands at unit speed in the normal direction while being advected by $V$, which models the wind velocity. When compared with homogenization of other Hamilton-Jacobi equations, the main difficulty with the G equation is that, since we do not assume that $\|V\|_{L^\infty} < 1$, the equation may not be coercive. On the other hand, if $\E[V] = 0$, then the equation is still ``coercive on average'', so we can recover some large-scale controllability.

Cardaliaguet\textendash{}Souganidis proved, under the assumption that the environment $V$ is stationary ergodic, that the equation homogenizes; i.e.\ we have the locally uniform convergence of solutions $u^\varepsilon \to \overline{u}$ as $\varepsilon \to 0$ almost surely, where $\overline{u}$ is the solution to the effective equation
\begin{equation}\label{macro-g-eqn}
    \begin{cases}
        D_t \overline{u}(t, x) = \overline{H}(D_x\overline{u}(t, x)) &\quad \text{ in $\R_{>0} \times \R^d$}\\
        \overline{u}(0, x) = u_0(x) &\quad \text{ in $\R^d$},
    \end{cases}
\end{equation}
and $\overline{H} \colon \R^d \to \R^d$ is a degree $1$ positively homogeneous coercive function, called the effective Hamiltonian.

Under the purely qualitative ergodicity assumption, there is no hope of proving a rate at which $u^\varepsilon$ converges to $\overline{u}$. In this paper, we make the stronger assumption that $V$ has unit range of dependence, which is a continuous analogue of i.i.d. Using this more quantitative assumption, we classify regions as ``good'' if the controllability estimate from Cardaliaguet\textendash{}Souganidis holds locally, in a quantitative sense. Then, we use percolation estimates to construct paths which stay inside the good regions. Our main result is the following rate of homogenization.
\begin{thm}\label{main-homog-theorem}
    Let $\P$ be a probability measure on $C^{1,1}(\R^d; \R^d)$ which has unit range of dependence and is $\R^d$-translation invariant. Assume that $\div V = 0$ and $\|V\|_{C^{1,1}} \leq L$ almost surely for some $L > 0$. Then there are constants $C(d, L) > 1 > c(d, L) > 0$ and a random variable $T_0$, with $\E[\exp(c\log^{3/2} T_0)] \leq C$, such that \[ |u^{\varepsilon}(t, x) - \overline{u}(t, x)| \leq C\|u_0\|_{C^{0,1}}{(t\varepsilon)}^{1/2}\log^3(\varepsilon^{-1}t) \] for all $t \geq \varepsilon T_0$ and $|x| \leq t$.
\end{thm}
Note that, in particular, the bound on $T_0$ implies $\E[T_0^n] \leq C^{2^n}$ for all $n \in \N$. In the interest of completeness, we record the almost-sure version of the rate of homogenization, which follows immediately from Theorem~\ref{main-homog-theorem} and the Borel\textendash{}Cantelli lemma.
\begin{cor}
    Under the same assumptions as Theorem~\ref{main-homog-theorem}, there is a constant $C(d, L) > 0$ such that, for all $T > 0$, \[ \limsup_{\varepsilon \to 0^+} \sup_{(t, x) \in [0, T] \times B_T} \frac{|u^{\varepsilon}(t, x) - \overline{u}(t, x)|}{{(T\varepsilon)}^{1/2}\log^3(\varepsilon^{-1}T)} \leq C\|u_0\|_{C^{0,1}} \] almost surely.
\end{cor}

As an application of the quantitative rate, we prove that the effective Hamiltonian depends continuously on the law of the environment.
\begin{thm}\label{continuous-dependence-theorem}
    Let $\P, {\{\P^n\}}_{n \in \N}$ be probability measures on $C^{1,1}(\R^d; \R^d)$ which are $\R^d$-translation invariant, have unit range of dependence, and satisfy $\div V = 0$ and $\|V\|_{C^{1,1}} \leq L$, $\P$-almost surely and $\P^n$-almost surely for every $n \in \N$. Let $\overline{H}^n$ and $\overline{H}$ be the effective Hamiltonians for $\P^n$ and $\P$ respectively. If $\P^n \weakstarto \P$ in the space of probability measures on $L^\infty(B_R; \R^d)$ for every $R > 0$, then $\overline{H}^n$ converges locally uniformly to $\overline{H}$.
\end{thm}

\subsection{Prior work} Since Lions\textendash{}Papanicolau\textendash{}Varadhan~\cite{LPV} proved qualitative homogenization of coercive Hamilton-Jacobi equations in a periodic environment, there has been a rich body of work (see Tran~\cite{Tran-book}) studying homogenization of Hamilton-Jacobi equations in both periodic and stochastic environments.

In a periodic environment, when the Hamiltonian is coercive, Capuzzo-Dolcetta\textendash{}Ishii~\cite{CDI} gave the first proof of a quantitative rate $O(\varepsilon^{1/3})$ of homogenization. Although the G equation may not be coercive when $\|V\|_{L^\infty} \geq 1$, it has a particularly simple optimal control formulation, making it an ideal first candidate to study homogenization of noncoercive equations. Cardaliaguet\textendash{}Nolen\textendash{}Souganidis~\cite{CardNoleSoug} proved homogenization, along with a quantitative rate $O(\varepsilon^{1/3})$, for the G equation in a periodic environment.

In a stochastic environment, the situation is more complicated. An example of Ziliotto~\cite{Ziliotto} shows that there exist coercive Hamiltonians and stationary ergodic environments in which homogenization does not hold. Feldman\textendash{}Souganidis~\cite{FeldSoug} generalized this example by showing that for any Hamiltonian with a strict saddle point, there exists a stationary ergodic environment in which homogenization does not hold. On the other hand, when the Hamiltonian is coercive and convex, Souganidis~\cite{Soug1999} proved that homogenization holds in a stationary ergodic environment. Later, under the stronger assumption that the environment has finite range of dependence, Armstrong\textendash{}Cardaliaguet\textendash{}Souganidis~\cite{ArmsCardSoug} proved a quantitative rate of homogenization when the Hamiltonian is coercive and level-set convex. In the case of the G equation, Cardaliaguet\textendash{}Souganidis~\cite{CardSoug} proved qualitative homogenization in a stationary ergodic environment.

In this paper, we combine ideas from~\cite{ArmsCardSoug} and~\cite{CardSoug} with percolation theory techniques to find a quantitative rate for the G equation, under the stronger finite range of dependence assumption.

\subsection{Assumptions} We now explicitly specify the assumptions in Theorems~\ref{main-homog-theorem} and~\ref{continuous-dependence-theorem}. Let $\P$ be a probability measure over $\Omega := C^{1,1}(\R^d; \R^d)$ with $d \geq 2$. We write $V \colon \Omega \to \Omega$ to denote the identity random variable, so $V$ has distribution $\P$. We assume that, $\P$-almost surely, $V$ is divergence-free and ${\|V\|}_{C^{1,1}} \leq L$, where $L > 0$ is a deterministic constant. Given $A \subseteq \R^d$, we write $\mathcal{G}(A)$ to denote the $\sigma$-algebra generated by $V$ restricted to $A$. That is, $\mathcal{G}(A)$ is the smallest $\sigma$-algebra for which the random variables $V(x)$ are measurable for every $x \in A$. We assume that $\P$ has unit range of dependence, which means that if $A, B \subseteq \R^d$ with $\dist(A, B) > 1$, then the $\sigma$-algebras $\mathcal{G}(A)$ and $\mathcal{G}(B)$ are independent. We also assume that $\P$ is $\R^d$-translation invariant, which means that $V(\cdot+x)$ has the same distribution as $V$ for any $x \in \R^d$. Note that this assumption is not much different than $\Z^d$-translation invariance for our purposes; we can simply add a random vector in ${[0, 1]}^d$ to go from $\Z^d$-invariance to full $\R^d$-invariance, which will not have an effect on our results. Finally, we assume that $\E[V] = 0$ (this is the same as $\E[V(0)] = 0$ by $\R^d$-translation invariance).

\subsection{Structure of the paper} In Section~\ref{percolation-section}, we collect some estimates from percolation theory. The bulk of the paper is Section~\ref{metric-problem-section}, in which we prove homogenization of the metric problem, which we later apply via the optimal control formulation. Section~\ref{metric-problem-section} has three main parts. First, we prove a controllability estimate for the metric problem, using results of Cardaliaguet\textendash{}Souganidis~\cite{CardSoug} and percolation estimates from Section~\ref{percolation-section}. With the controllability estimate in hand, we split the difference between the microscopic and macroscopic metric problems into two pieces; the random fluctuations and the nonrandom scaling bias. We handle the first using a martingale argument originally due to Kesten~\cite{Kesten93}, and the second with an argument adapted from Alexander~\cite{Alexander97}. Finally, we apply our estimates for the metric problem in Section~\ref{applications-section} to deduce our main results.

\subsection{Notation} Throughout the paper, $C > 1 > c > 0$ will denote constants which may depend on the dimension $d$ and the bound $L$ on $\|V\|_{C^{1,1}}$, but may vary from line to line. We write $B_r := \{x \in \R^d \mid |x| \leq r\}$ for the Euclidean ball of radius $r$ centered at the origin. We write $\cl(\cdot)$ (resp. $\intr(\cdot)$) to denote the closure (resp.\ interior) of a subset of some topological space. The Euclidean distance between two subsets $E, F \subseteq \R^d$ is given by $\dist(E, F) := \inf \{|x-y| \mid (x, y) \in U \times V\}$. The Hausdorff distance between two subsets $E, F \subseteq \R^d$ is given by \[ \dist_H(E, F) := \inf \{\varepsilon > 0 \mid V \subseteq U + B_\varepsilon \text{ and } U \subseteq V + B_\varepsilon\}. \] For a random variable $X$, we use the subscript $X_\omega$ to denote the value of $X$ at $\omega \in \Omega$.

\subsection{Acknowledgements} I would like to thank my advisor, Charles Smart, for suggesting the problem and many helpful conversations. I would also like to thank Panagiotis Souganidis for suggesting the problem of continuous dependence of the Hamiltonian on the law of the environment.

\numberwithin{thm}{section}
\numberwithin{cor}{section}

\section{Percolation estimates}\label{percolation-section}
In this section, we prove a few well-known estimates from supercritical percolation theory, but with a finite range of dependence assumption. This context is similar to the usual one, where the environment is i.i.d., except that the underlying probability space lacks a product structure.

Let $d \geq 2$ and let $G \colon \Z^d \to \{0, 1\}$ be a random function on $\Z^d$ which has finite range $C_{\text{dep}} > 0$ of dependence, which means that the $\sigma$-algebras induced by the values of $G$ on sets that are Euclidean distance at least $C_{\text{dep}}$ apart are independent. We assume that $G$ is $\Z^d$-translation invariant, i.e., $G(\cdot+v)$ has the same distribution as $G(\cdot)$ for all $v \in \Z^d$. The function $G$ models site percolation, where a site $x$ is open if $G(x) = 1$ and closed otherwise.

We put edges on $\Z^d$ between nearest neighbors in the $\ell^\infty$ metric. We write $\dist(\cdot, \cdot)$ to indicate the graph distance, and refer to maximal connected components on which $G$ is constant as \textit{clusters}. Clusters composed of open (resp.\ closed) sites are called open (resp.\ closed) clusters. Let $p := \P[G(0) = 1]$ be the probability that a site is open (which is the same for every site, by $\Z^d$-translation invariance). We write $\mathcal{Q}_R(x) \subseteq \Z^d$ to denote the axis-aligned cube of side length $2R$ centered at $x$.
\begin{lem}\label{closed-clusters-small}
    Let $S \subseteq \Z^d$ be a finite set. Define the closed sites connected to $S$ by \[ \cl(S) := \left\{x \in \Z^d \mid \text{there is a path of closed sites from $x$ to a site in $S$}\right\}. \] For any $\varepsilon > 0$ there is $p_0(d, C_\text{dep}, \varepsilon) < 1$ and $C(d, C_\text{dep}, \varepsilon) > 0$ such that if $p > p_0$ then \[ \P\left[|\cl(S)| > \varepsilon|S| + \delta\right] \leq C\exp(-C^{-1}\delta). \]
\end{lem}
Note in particular that $x \in \cl(S)$ implies that $x$ is a closed site.
\begin{proof}
    Let $T \subseteq \Z^d$ be a finite set of $n$ vertices. If $T \subseteq \cl(S)$, then every site in $T$ is closed and every cluster in $T$ contains a point in $S$. For fixed $n$, the number of sets $T$ which satisfy the latter condition is at most ${(2d+2)}^{|S|+2n}$ (we can encode a spanning tree of $T$ with an alphabet of $2d+2$ letters). For a fixed set $T$, we see that \[ \P[\text{every site in $T$ is closed}] \leq {(1-p)}^{\left\lfloor n/{(2C_\text{dep})}^d\right\rfloor}, \] since we can choose at least $\left\lfloor n/{(2C_\text{dep})}^d\right\rfloor$ sites in $T$ which are far enough to be independent. From the union bound, we have \[\P[\text{there is a set $T \subseteq \cl(S)$ with $n$ vertices}] \leq {(2d+2)}^{|S|+2n}{(1-p)}^{\left\lfloor n/{(2C_\text{dep})}^d\right\rfloor}. \] Now let $n := \left\lceil\varepsilon|S| + \delta\right\rceil$ and choose $(1-p_0)$ small enough so that \[{(1-p_0)}^{\varepsilon/{(2C_\text{dep})}^d}{(2d+2)}^{1+2\varepsilon} < 1\] and \[{(1-p_0)}^{1/{(2C_\text{dep})}^d}{(2d+2)}^2 < 1\] and the claim follows.
\end{proof}
The next lemma has nothing to do with the percolation environment; it's simply a property of the graph structure of $\Z^d$. It follows from a topological property of $\R^d$ known as unicoherence (see Kuratowski~\cite{Kuratowski} or Dugundji~\cite{Dugundji}). In order to state the lemma, we need to define the boundary of a subset of $E \subseteq \Z^d$. Because $\Z^d$ is discrete, there are two choices for our definition.
\begin{defn}
    The inner (resp.\ outer) boundary of $E$, denoted $\partial^-E$ (resp. $\partial^+E$), is the set \[ \partial^-E := \{x \in E \mid \dist(x, \Z^d \setminus E) = 1\} \qquad \text{(resp.\ $\partial^+E := \{x \in \Z^d \setminus E \mid \dist(x, E) = 1\}$)}. \]
\end{defn}

\begin{lem}\label{unicoherence}
    Let $\mathcal{Q}_R$ be any cube of side length $2R$ and let $\mathfrak{C} \subseteq \mathcal{Q}_R$ be a connected set. Let $\mathfrak{D} \subseteq \mathcal{Q}_R \setminus \mathfrak{C}$ be a connected component of $\mathcal{Q}_R \setminus \mathfrak{C}$. Then the inner (resp.\ outer) boundary of $\mathfrak{D}$ is connected.
\end{lem}
\begin{proof}
This is part (i) of Lemma 2.1 from Deuschel\textendash{}Pisztora~\cite{DeusPisz}. The proof is a standard application of Urysohn's lemma.
\end{proof}

The next lemma shows that, with high probability, there is a large open cluster which is near every site.
\begin{lem}\label{open-cluster-big}
    Let $n, R > 0$ and consider $\mathcal{Q}_R$, a cube of side length $2R$. Let $E_n$ be the event that there exists an open cluster $\mathfrak{C} \subseteq \mathcal{Q}_{R+n}$ such that every connected component of $\mathcal{Q}_{R+n} \setminus \mathfrak{C}$ which intersects $\mathcal{Q}_R$ is of size at most $n$. Then there are constants $p_0(d, C_\text{dep}) < 1$ and $C = C(d, C_\text{dep}) > 0$ such that if $p > p_0$ then \[\P[E_n] \geq 1-CR^d\exp(-C^{-1}n^{(d-1)/d}).\]
\end{lem}
\begin{proof}
    Fix $n, R > 0$ as in the statement. Work in the event that every closed cluster in $\mathcal{Q}_{R+n}$ has size less than $C^{-1}n^{(d-1)/d}$, where $C = C(d)$ comes from the isoperimetric constant (chosen later in the proof). By Lemma~\ref{closed-clusters-small} (applied to each site in $\mathcal{Q}_{R+n}$ individually, with (say) $\varepsilon = 1$), this event has probability at least $1-CR^d\exp(-C^{-1} n^{(d-1)/d})$.

    Let $\mathfrak{C}$ be the largest open cluster (breaking ties arbitrarily) in $\mathcal{Q}_{R+n}$. As long as $C \geq 1$, it follows from the isoperimetric inequality that there is an open path between opposite faces of $\mathcal{Q}_{R+n}$, so $|\mathfrak{C}| \geq 2(R+n)+1 > n$.

    Let $\mathfrak{D}$ be any connected component of $\mathcal{Q}_{R+n} \setminus \mathfrak{C}$ which intersects $\mathcal{Q}_R$. The inner boundary of $\mathfrak{D}$ is composed of two kinds of sites: (i) those bordering $\mathfrak{C}$ and (ii) those in the inner boundary of $\mathcal{Q}_{R+n}$. The sites of type (i) are all closed (by definition of $\mathfrak{C}$). We claim that there are no sites of type (ii). Indeed, if there was a site of type (ii) then we could follow the inner boundary of $\mathfrak{D}$ (it's connected by Lemma~\ref{unicoherence}) from the inner boundary of $\mathcal{Q}_{R+n}$ to a site in $\mathcal{Q}_R$, which would yield a path of length more than $n \geq n^{(d-1)/d}$ of type (i) (and hence closed) sites, contradicting our assumption.

    Since the inner boundary of $\mathfrak{D}$ is connected and composed entirely of closed sites, it has size less than $C^{-1}n^{(d-1)/d}$. The isoperimetric inequality then shows that either $\mathfrak{D}$ or $\mathcal{Q}_{R+n} \setminus \mathfrak{D}$ has size at most $n$. Since $\mathfrak{C} \subseteq \mathcal{Q}_{R+n} \setminus \mathfrak{D}$ and $|\mathfrak{C}| > n$, it follows that $|\mathfrak{D}| \leq n$ as desired.
\end{proof}

\section{The metric problem}\label{metric-problem-section}
We now shift our focus to the metric problem associated with the G equation, which comes from the optimal control formulation. As we will see, homogenization of solutions to the G equation is implied by convergence of the associated metric to its large-scale limit.

Given $t > 0$ and a measurable function $\alpha \colon [0, t] \to B_1$, define the \textit{controlled path} $X_x^\alpha \colon [0, t] \to \R^d$ to be the solution to the initial-value problem
\begin{equation}\label{controlled-path-ode}
    \begin{cases} \dot{X}_x^\alpha = \alpha + V(X_x^\alpha)\\ X_x^\alpha(0) = x. \end{cases}
\end{equation}
For each $x \in \R^d$, define the \textit{reachable set at time $t$} by
\begin{equation}\label{reachable-set-defn}
    \mathcal{R}_t(x) := \{y \in \R^d \mid \exists\; \alpha \colon [0, t] \to B_1 \text{ such that } X_x^\alpha(t) = y\}.
\end{equation}
Note that this definition still makes sense for $t < 0$, if we interpret $[0, t]$ as $[t, 0]$. For convenience, we also define the sets \[ \mathcal{R}_t^-(x) := \bigcup_{0 \leq s \leq t} \mathcal{R}_s(x) \] for $t \geq 0$ and \[ \mathcal{R}_t^+(x) := \bigcup_{t \leq s \leq 0} \mathcal{R}_s(x) \] for $t \leq 0$. Define the \textit{first passage time}
\begin{equation}\label{theta-definition}
    \theta(x, y) := \inf \{t \mid y \in \mathcal{R}_t(x)\}.
\end{equation}
Finally, if $E \subseteq \R^d$ is a set, we define \[ \mathcal{R}_t(E) = \bigcup_{e \in E} \mathcal{R}_t(e), \] and we do the same for $\mathcal{R}_t^-$ and $\mathcal{R}_t^+$.

\subsection{Controllability}

First, we import a few results from Cardaliaguet\textendash{}Souganidis~\cite{CardSoug}, which hold in the more general ergodic setting.
\begin{lem}[Cardaliaguet\textendash{}Souganidis~\cite{CardSoug}, Lemma 4.2]\label{CS-growth}
    There is a deterministic constant $\beta  = \beta(d) > 0$ such that $\mathcal{R}_t(x) \geq \beta|t|^d$ (and hence the same holds for $\mathcal{R}_t^-(x)$ and $\mathcal{R}_t^+(x)$).
\end{lem}
\begin{thm}[Cardaliaguet\textendash{}Souganidis~\cite{CardSoug}, Theorem 4.1]\label{CS-control}
    For every $\varepsilon > 0$, there is a random variable $0 < T(\varepsilon) < \infty$ such that \[ \theta(x, y) \leq T(\varepsilon) + \varepsilon|x| + (1+\varepsilon)|x-y| \qquad \text{for all $x, y \in \R^d$}\] holds almost surely.
\end{thm}

For our quantitative purposes, we would like to have a version of Theorem~\ref{CS-control} where the distribution of the error term $T(\varepsilon)$ has exponential tail bounds. Our strategy is to partition space into cubes where a controllability estimate holds nearby with high probability. Then, using the percolation estimates, we can construct paths which mostly stay in these cubes.
\begin{lem}\label{waiting-time-exists}
    For each $0 < p < 1$ there is $C = C(d, p, L) \in \N$ such that \[ \P\left[\sup_{x, y \in B_{\sqrt{d}}} \theta(x, y) \leq C\right] \geq p. \]
\end{lem}
\begin{proof}
    By Theorem~\ref{CS-control}, the statement of the lemma holds when $C$ is also allowed to depend on the distribution of $V$, since $T(\varepsilon)$ is almost surely finite. We claim that we can remove this dependence. Indeed, if not, then there would be some sequence ${\{\P_n\}}_{n \in \N}$ of probability measures which satisfy the same assumptions as $\P$, as well as \[ \P_n\left[E(n)\right] < p, \] where $E(C)$ is the event that $\theta(x, y) \leq C$ for all $x, y \in \overline{B_{\sqrt{d}}}$. Write $\tilde{\theta}(x, y)$ to denote the first-passage time where the control $\alpha$ is constrained to lie in $B_{1/2}$ instead of $B_1$ as in the definition~(\ref{reachable-set-defn}) of the reachable set. Similarly, define $\tilde{E}(C)$ to be the event that $\tilde{\theta}(x, y) \leq C$ for all $x, y \in \overline{B_{\sqrt{d}}}$, noting that \[ \cl(\tilde{E}(C)) \subseteq \intr(E(C)), \] where $\cl(\cdot)$ and $\intr(\cdot)$ denote the closure and interior respectively, taken in the space $C^1(\R^d; \R^d)$. Taking a subsequence (not relabelled for brevity), we find a probability measure $\P_\infty$, which satisfies all the same assumptions as $\P$, such that $\P_n|_{B_R} \weakstarto \P_\infty|_{B_R}$ in the space of probability measures on $C^1(B_R; \R^d)$ for every $R > 0$. Then, for any $C > 0$, \[ p > \limsup_{n \to \infty} \P_n[E(C)] \geq \P_\infty[\tilde{E}(C)], \] which violates Theorem~\ref{CS-control} (noting that Theorem~\ref{CS-control} holds just as well for $\tilde{\theta}$ as for $\theta$ by a change of variables).
\end{proof}

\begin{lem}\label{it-is-percolation}
    For each $0 < p < 1$, there exists a constant $C = C(d, p, L) > 0$ such that the function $G \colon \Z^d \to \{0, 1\}$, defined by \[ G(v) = \begin{cases} 1 & \quad \theta(x, y) \leq C \; \text{for all $x, y \in B_{\sqrt{d}}(v)$}\\ 0 & \quad \text{otherwise,} \end{cases} \] is $Z^d$-translation invariant with finite range of dependence $C_\text{dep} = C_\text{dep}(d, p, L) > 0$ and $\P[G(0) = 1] \geq p$.
\end{lem}
 In other words, $G$ is an environment in which all of our percolation estimates apply.
\begin{proof}
    By Lemma~\ref{waiting-time-exists}, there is $C > 0$ such that $\P[G(0) = 1] \geq p$. The $\Z^d$-translation invariance of $G$ follows from that of $\P$. Finite range of dependence follows from the fact that the value of $G(v)$ depends only on controlled paths starting in $B_{\sqrt{d}}(v)$ which run for time at most $C$. Since $\|V\|_{L^\infty} \leq L$, these paths stay inside of $\overline{B_{\sqrt{d} + (L+1)C}(v)}$, so $C_\text{dep} := 2(\sqrt{d} + (L+1)C) + 1$ satisfies the claim.
\end{proof}

\begin{defn}
    To translate between $\Z^d$ and $\R^d$, for each set $E \subseteq \Z^d$ we introduce the ``solidification'' \[\sigma(E) := E+{\left[-\frac12, \frac12\right]}^d. \]
\end{defn}

 We are now in a position to prove our improved controllability estimate. We refer to clusters as before, using the same notion of adjacency. We write $W_R := \sigma(\mathcal{Q}_R)$
\begin{thm}\label{cont-estimate}
    For each $R \geq 1$, there is a constant $C = C(d, L) > 0$ such that the ``extra waiting time'' \[ \mathcal{E}(R) := \sup_{x, y \in W_{R}} \theta(x, y) - (C + C|x-y|) \] satisfies \[ \P[\mathcal{E}(R) > n] \leq CR^d\exp(-C^{-1} n). \]
\end{thm}
\begin{proof}
    We partition $\R^d$ into cubes of side length $1$, centered at points in $\Z^d$. Given Lemma~\ref{it-is-percolation}, there is a constant $C = C(d, p, L) > 0$ such that $\P[G(v) = 1] \geq p$ for some $p > p_0(d, C_\text{dep}, \varepsilon)$, where $\varepsilon = 1$ will be chosen later. For $v \in \Z^d$, we say that the site $v$ is open if $G(v) = 1$ and closed otherwise. We say that a point $x \in \R^d$ lies near an open site if there is some $v \in \Z^d$ such that $x \in \sigma(\{v\})$. Let $S = \lceil R \rceil$ and $x, y \in W_{S}$. For convenience, we will prove that $\P[\mathcal{E}(R) > Cn] \leq CR^d\exp(-C^{-1}n)$; this easily implies the original claim by changing $C$ to $C^2$.

    \textit{Step 1.} We show that we can assume without loss of generality that $x$ and $y$ lie in the same open cluster, by which we mean that there is an open cluster $\mathfrak{C}$ such that $x, y \in \sigma(\mathfrak{C})$. Indeed, suppose they lie in different open clusters. Then by Lemma~\ref{open-cluster-big}, we have an open cluster $\mathfrak{C} \subseteq \mathcal{Q}_{S+n^d}$, which depends on the environment, such that \[ \P\left[\text{every connected $\mathfrak{D} \subseteq \mathcal{Q}_{S+n^d} \setminus \mathfrak{C}$ which intersects $\mathcal{Q}_S$ satisfies $|\mathfrak{D}| \leq n^d$}\right] \geq 1-CS^d\exp(-C^{-1}n^{d-1}).\] Working in this event, let $\mathfrak{D}$ be the connected component of $\mathcal{Q}_{S+n^d} \setminus \mathfrak{C}$ whose solidification $\sigma(\mathfrak{D})$ contains $x$. Then Lemma~\ref{CS-growth}, along with the fact that $\mathcal{R}_t^-(x)$ is connected, shows that $\mathcal{R}_t^-(x) \setminus \sigma(\mathfrak{D}) \neq \emptyset$ if $\beta|t|^d \geq n^d$ and $t > 0$. Since the controlled paths are continuous, there is some $t \leq n\beta^{-1/d} + 1$ such that there is some $z \in \mathcal{R}_t^-(x) \cap \sigma(\mathfrak{C})$. Since the loss in probability and travel time can be controlled by enlarging $C$, we may as well assume that $x$ was $z$ to begin with. Repeating the same argument for $y$ (except running time backwards, so $t < 0$ and we use $\mathcal{R}_t^+(y)$ instead) shows that we may assume that $x$ and $y$ lie in the same open cluster, $\mathfrak{C}$.

    \textit{Step 2.} To show that $\theta(x, y) - C|x-y| \leq Cn$, we will build a ``skeleton'' of points $x = x_0, x_1, \dots, x_k = y$ which all lie near open sites and satisfy $|x_{i+1}-x_i| \leq \sqrt{d}$ for each $0 \leq i < k$. Then, by connecting the points with paths given by Lemma~\ref{it-is-percolation}, we can build a controlled path of length at most $Ck$ which follows the skeleton. It remains to show that we can build such a skeleton with $k \leq C(|x-y| + n)$. Let \[A := \{v \in \Z^d \mid \sigma(\{v\}) \cap \overline{xy} \neq \emptyset\}\] be the set of centers of cubes which intersect the line segment connecting $x$ and $y$. Note that $|A| \leq 2^d(1+|x-y|)$. By Lemma~\ref{closed-clusters-small}, we can choose $\varepsilon = 1$ to work in the event that $\cl(A) \leq |A| + n$.

    Our strategy is to go from $x$ to $y$ in a straight line, taking necessary detours around closed clusters. We use Lemma~\ref{closed-clusters-small} to bound the total length of our detour.

    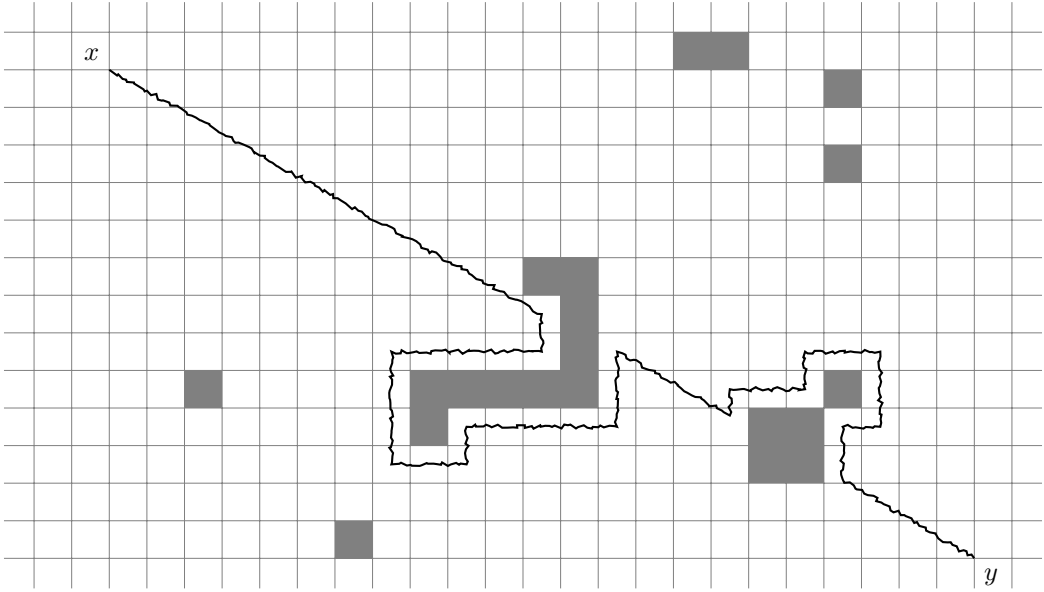
\begin{figure}
        \centering
        \begin{tikzpicture}
            \usetikzlibrary{decorations.pathmorphing}
            \draw[step=0.5cm, gray, very thin] (-6.9, -3.9) grid (6.9, 3.9);
            \fill[gray] (0, 0) rectangle (0.5, 0.5);
            \fill[gray] (0.5, 0) rectangle (1, 0.5);
            \fill[gray] (2, 3) rectangle (2.5, 3.5);
            \fill[gray] (2.5, 3) rectangle (3, 3.5);
            \fill[gray] (-1.5, -2) rectangle (-1, -1.5);
            \fill[gray] (0.5, -0.5) rectangle (1, 0);
            \fill[gray] (0.5, -1) rectangle (1, -0.5);
            \fill[gray] (0.5, -1.5) rectangle (1, -1);
            \fill[gray] (0, -1.5) rectangle (0.5, -1);
            \fill[gray] (-0.5, -1.5) rectangle (0, -1);
            \draw[black, thick, decorate, decoration={random steps, segment length=0.05cm, amplitude=0.03cm}] (-5.5, 3) node[anchor=south east] {$x$} -- (0.25, -0.25) -- (0.25, -0.75) -- (-1.75, -0.75) -- (-1.75, -2.25) -- (-0.75, -2.25) -- (-0.75, -1.75) -- (1.25, -1.75) -- (1.25, -0.75) -- (2.75, -1.6) -- (2.75, -1.25) -- (3.75, -1.25) -- (3.75, -0.75) -- (4.75, -0.75) -- (4.75, -1.75) -- (4.25, -1.75) -- (4.25, -2.45) -- (6, -3.5) node[anchor=north west] {$y$};
            \fill[gray] (-1, -1.5) rectangle (-0.5, -1);
            \fill[gray] (-1.5, -1.5) rectangle (-1, -1);
            \fill[gray] (-2.5, -3.5) rectangle (-2, -3);
            \fill[gray] (-4.5, -1.5) rectangle (-4, -1);
            \fill[gray] (4, -1.5) rectangle (4.5, -1);
            \fill[gray] (3.5, -2) rectangle (4, -1.5);
            \fill[gray] (3, -2) rectangle (3.5, -1.5);
            \fill[gray] (3.5, -2.5) rectangle (4, -2);
            \fill[gray] (3, -2.5) rectangle (3.5, -2);
            \fill[gray] (4, 2.5) rectangle (4.5, 3);
            \fill[gray] (4, 1.5) rectangle (4.5, 2);
        \end{tikzpicture}
        \caption{An example of our controlled path from $x$ to $y$; cubes corresponding to closed sites are shaded}
    \end{figure}

    We can build the skeleton iteratively. Start with $x_0 = x$, and assume we've built the skeleton up to $x_i$, for $i \geq 0$. We maintain the invariant that, at the end of each step, $x_i$ lies on the line segment $\overline{xy}$, and it is closer to $y$ than any other $x_j$ which lies on $\overline{xy}$ with $0 \leq j < i$. There are three cases.
    \begin{enumerate}
        \item If $|y-x_i| \leq \sqrt{d}$, define $x_{i+1} := y$ and finish.
        \item If $z := x_i + \sqrt{d}\frac{y-x}{|y-x|}$ lies near an open site, define $x_{i+1} := z$ and continue to the next step.
        \item Otherwise, let $\tilde{x}$ be the point on $\overline{x_i z}$ which lies near an open site and is closest to $z$. Then $\tilde{x}$ also lies near a closed site, which is part of some connected component $\mathfrak{F}$ of $\mathcal{Q}_S \setminus \mathfrak{C}$. By Lemma~\ref{unicoherence}, the outer boundary of $\mathfrak{F}$ is connected. Since $\tilde{x} \in \overline{xy}$, we see that $\partial^- \mathfrak{F} \subseteq \cl(A)$. Besides, since every vertex in $\Z^d$ has degree $3^d-1$, we have $|\partial^+ \mathfrak{F}| \leq 3^d |\partial^- \mathfrak{F}|$. So, let $p_1, \dots, p_\ell \in \Z^d$ be a path along the outer boundary of $\mathfrak{F}$, where $\sigma(\{p_1\}) \ni \tilde{x}$ and $p_\ell \in A$ is a point on the outer boundary of $\mathfrak{F}$ which maximizes $p_\ell \cdot (y-x)$. Finally, we extend our path by setting $x_{i+1} := \tilde{x}$ and $x_{i+j+1} := p_j$ for each $1 \leq j \leq \ell$, and set $x_{i+\ell+2}$ to be a point on the line segment $\overline{xy}$ which lies in $\sigma(\{p_\ell\})$.
    \end{enumerate}

    It remains to analyze the length of this skeleton by looking at each of the three cases. The first case happens at most once, so it can be ignored. The second case reduces the distance $|x_i-y|$ by $\sqrt{d}$ and the third case does not increase the distance $|x_i-y|$, so there can be at most $\frac{|x-y|}{\sqrt{d}}$ points in the skeleton coming from the second case. The third case adds $\ell+2$ points, where $\ell \leq |\partial^+\mathfrak{F}|$. Since we finish an instance of the third case at a point as close to $y$ as possible on the segment $\overline{xy}$, we never witness the same cluster $\mathfrak{F}$ twice in different instances of the third case. Therefore the third case adds at most \[ C|\cl(A)| \leq C(|A| + n) \leq C(|x-y| + n) \] points to our skeleton.
\end{proof}

\subsection{Random fluctuations in first passage time}

Next, we consider how much $\theta(0, y)$ deviates from its expectation. Our proof will follow roughly the same path as the proof of Proposition 4.1 from Armstrong\textendash{}Cardaliaguet\textendash{}Souganidis~\cite{ArmsCardSoug}, with some modifications which are made possible by the controllability estimate.

To get started, we introduce a ``guaranteed'' version of first passage time. For any $\rho > 0$, define the $\rho$-guaranteed reachable set recursively by \[ \mathcal{R}_t^\rho(x) := \begin{cases} \mathcal{R}_t^-(x) & \quad \text{if $t < \rho$}\\ \mathcal{R}_\rho^-(\mathcal{R}_{t-\rho}^\rho(x)) \cup (\mathcal{R}_{t-\rho}^\rho(x) + \overline{B_1}) & \quad \text{otherwise.} \end{cases} \] The $\rho$-guaranteed reachable set is similar to the reachable set, except that we enforce expansion at a rate of at least $1/\rho$ in a certain discrete sense. We similarly define the $\rho$-guaranteed first passage time \[ \theta^\rho(x, y) = \min \{t \geq 0 \mid y \in \mathcal{R}_t^\rho(x) \}. \]

Note that the $\rho$-guaranteed first passage time coincides with the usual first passage time if we have sufficient control on the extra waiting time $\mathcal{E}$ (from Theorem~\ref{cont-estimate}) in a suitable domain.

Fix some $y \in \R^d$ and define the random variable ${\{Z^\rho_t\}}_{t \geq 0}$ by \[ Z^\rho_t := \E\left[\theta^\rho(0, y) \mid \mathcal{F}_t\right], \] where $\mathcal{F}_t$ is the $\sigma$-algebra generated by the environment $V(x)$ restricted to the $\rho$-guaranteed reachable set $\mathcal{R}_t^\rho(0)$. In other words, $\mathcal{F}_t$ is the smallest $\sigma$-algebra so that the functions $V(x)\mathds{1}_{x \in \mathcal{R}_t^\rho(0)}$ are $\mathcal{F}_t$-measurable for every $x \in \R^d$. Since $\mathcal{R}_t^\rho(0)$ are increasing sets, ${\{\mathcal{F}_t\}}_{t \geq 0}$ is a filtration, so ${\{Z^\rho_t\}}_{t \geq 0}$ is a martingale.

We first show that $Z^\rho_t$ depends mostly on the shape of $\mathcal{R}_t^\rho(0)$, without regard for the values of $V$ inside $\mathcal{R}_t^\rho(0)$. In order to condition on the approximate shape of the reachable set, for any $E \subseteq \R^d$ we introduce the discretization \[ \disc(E) := \{z \in d^{-1/2}\Z^d \mid B(z, 1) \cap E \neq \emptyset\}. \]
\begin{lem}\label{fluc-helper-1}
    For any $t \geq 0$, we have \[ \left|\max(Z^\rho_t, t) - f(t, \disc(\mathcal{R}_t^\rho(0)))\right| \leq 3\rho, \] where we define $f(t, S)$, for any $t \geq 0$ and any finite set $S \subseteq d^{-1/2}\Z^d$, by \[ f(t, S) := t + \E\left[\theta^\rho(S, y)\right]. \]
\end{lem}
\begin{proof}
    Fix some $t \geq 0$. Using the speed limit $L+1$ for controlled paths, we see that \[ \mathcal{R}_t^\rho(0) \subseteq B(0, 1+\lceil L + 1 + \rho^{-1} \rceil t) \] almost surely. Define the set of possible discretized reachable sets at time $t$ by
    \begin{equation}\label{C_t-definition}
        C_t := \left\{S \subseteq d^{-1/2}\Z^d \cap B\left(0, \lceil L + 1 + \rho^{-1} \rceil t\right)\right\}.
    \end{equation}
    For any $S \in C_t$, the event that $\disc(\mathcal{R}^\rho_t(0)) = S$ is $\mathcal{F}_t$-measurable, so we have \[ \left|\max\left(Z^\rho_t, t\right) - f\left(t, \disc(\mathcal{R}_t^\rho(0))\right)\right| = \sum_{S \in C_t} \max\left(\left|\E\left[\left(\theta^\rho(0, y)-f(t, S)\right)\mathds{1}_{\disc(\mathcal{R}^\rho_t(0)) = S} \mid \mathcal{F}_t\right]\right|, \left|t-f(t, S)\right|\right). \]
    Fix some $S \in C_t$. If $B(y, 2) \cap S \neq \emptyset$, then $y \in \mathcal{R}^\rho_{2\rho}(S)$, so $\theta^\rho(0, y) \leq t+3\rho$ and $t \leq f(t, S) \leq t+3\rho$ and the conclusion holds.

    Otherwise, define the set $E := S + B(0, 2)$. Note that $\mathcal{R}^\rho_t(0) \subseteq E$ and $\dist(\mathcal{R}_t^\rho(0), \partial E) \geq 1$. Using the definition of the $\rho$-guaranteed reachable set,
    \begin{equation}\label{theta-p-inequality}
        \theta^\rho(0, \partial E) + \theta^\rho(\partial E, y) - \rho \leq \theta^\rho(0, y) \leq \theta^\rho(0, \partial E) + \theta^\rho(\partial E, y).
    \end{equation}
    Using the definitions of $S$ and $E$, \[ t \leq \theta^\rho(0, \partial E) \leq t + 3\rho. \] On the other and, the term $\theta^\rho(\partial E, y)$ is $\mathcal{G}(\R^d \setminus E)$-measurable. Taking the conditional expectation of~(\ref{theta-p-inequality}), we have \[ t - \rho + \E[\theta^\rho(\partial E, y)] \leq Z^\rho_t \leq t + 3\rho + \E[\theta(\partial E, y)]. \] To finish, we use the definition of the $\rho$-guaranteed reachable set to find that \[ 0 \leq \theta^\rho(S, y) - \theta^\rho(\partial E, y) \leq 2\rho. \] Combining the previous two displays yields the conclusion of the lemma.
\end{proof}

Next, we show that our approximation for $Z_t^\rho$, given by $f(t, \disc(\mathcal{R}_t^\rho(0)))$, has bounded increments.
\begin{lem}\label{fluc-helper-2}
    Let $t, s \geq 0$. Then \[ \left|f(t, \disc(\mathcal{R}_t^\rho(0))) - f(s, \disc(\mathcal{R}_s^\rho(0)))\right| \leq 2\rho + |t-s|(L\rho + \rho + 2). \]
\end{lem}
\begin{proof}
    Without loss of generality, assume $s < t$. Let $C_s$ and $C_t$ be as defined in~(\ref{C_t-definition}). We need to prove that \[ \sum_{A_s \in C_s} \sum_{A_t \in C_t} \left|f(t, A_t) - f(s, A_s)\right|\mathds{1}_{\disc(\mathcal{R}_s^\rho(0)) = A_s}\mathds{1}_{\disc(\mathcal{R}_t^\rho(0)) = A_t} \leq 2\rho + |t-s|(L\rho + \rho + 2). \] So, fix any $A_s \in C_s$ and $A_t \in C_t$ such that \[ \P\left[\disc(\mathcal{R}_s^\rho(0)) = A_s \text{ and } \disc(\mathcal{R}_t^\rho(0)) = A_t\right] > 0.\] The speed limit for controlled paths shows that $A_s \subseteq A_t$ and \[ \dist_H(A_s, A_t) \leq 2 + |t-s|\lceil L+1+\rho^{-1}\rceil, \] where $\dist_H$ denotes the Hausdorff distance. Using the definition of the $\rho$-guaranteed reachable set, this yields \[ \theta^\rho(A_t, y) \leq \theta^\rho(A_s, y) \leq \theta^\rho(A_t, y) + \rho\left(2+|t-s|\lceil L+1+\rho^{-1}\rceil\right). \] The conclusion of the lemma follows from the definition of $f$.
\end{proof}

Now we put these lemmas together and apply Azuma's inequality to ${\{Z^\rho_t\}}_{t \geq 0}$, choosing $\rho$ carefully to balance competing error terms.
\begin{prop}\label{random-fluc}
    Let $y_1, y_2 \in \R^d$ and $\lambda \geq C + C{|y_1-y_2|}^{1/2}\log^2|y_1-y_2|$. Then there is a constant $C = C(d, L) > 0$ such that \[ \P[|\theta(y_1, y_2) - \E[\theta(y_1, y_2)]| > \lambda] \leq C\exp\left(\frac{-C^{-1} \lambda^{1/2}}{|y_1-y_2|^{1/4}}\right). \]
\end{prop}
\begin{proof}
    By $\R^d$-translation invariance, we assume without loss of generality that $y_1 = 0$; let $y := y_2$.

    Note that $\theta(0, y) = \theta^\rho(0, y)$ as long as $\theta(u, v) \leq \rho$ if $|u-v| \leq 1$, for all $|u|, |v| \leq \lceil (L+1+\rho^{-1}) \rceil \lceil \rho^{-1}|y| \rceil$ (a ball of this radius contains the $\rho$-guaranteed reachable set at time $\lceil \rho^{-1}|y|\rceil$). By Lemma~\ref{fluc-helper-1} and Lemma~\ref{fluc-helper-2}, the martingale ${\{Z^\rho_t\}}_{t \geq 0}$ has bounded increments of \[ |Z_t - Z_s| \leq (\rho L+\rho+2)|t-s| + 8\rho. \] Also, $Z_t^\rho = \theta(0,y)$ for all $t \geq \rho|y|$. Apply the union bound with Theorem~\ref{cont-estimate} and Azuma's inequality to the sequence ${\{Z^\rho_n\}}_{0 \leq n \leq \lceil \rho^{-1}|y| \rceil}$ to see that
    \begin{align*}
        \P\left[\left|\theta(0, y) - \E[\theta^\rho(0, y)]\right| > \lambda\right] &\leq \P[\theta(0, y) \neq \theta^\rho(0, y)] + \P\left[\left|\theta(0, y) - \E[\theta^\rho(0, y)]\right| > \lambda\right]\\
        &\leq C{\left[(L+1+\rho^{-1})\rho^{-1}|y|\right]}^d\exp(-C^{-1} \rho)\\
        &\qquad + 2\exp\left(\frac{-C^{-1}\lambda^2}{{\left[(\rho L+\rho+2) + 8\rho\right]}^2\left\lceil \rho^{-1}|y|\right\rceil}\right),
    \end{align*}
    as long as $\rho \geq 2C$.
    Choosing $\rho = \lambda^{1/2}|y|^{-1/4}$ yields the bound
    \begin{equation}\label{ok-but-replace}
        \P[|\theta(0, y) - \E[\theta^\rho(0, y)]| > \lambda] \leq C\exp\left(\frac{-C^{-1} \lambda^{1/2}}{|y|^{1/4}}\right),
    \end{equation}
    as long as $\lambda \geq C|y|^{1/2}\log^2|y|$ (we absorb all the polynomials into the exponential by changing the constant appropriately). It remains to replace $\E[\theta^\rho(0, y)]$ in~(\ref{ok-but-replace}) with $\E[\theta(0, y)]$. Indeed, we can bound
    \begin{align}
        \E[\theta(0, y)] &= \E\left[\theta^\rho(0, y) \mid E_\rho\right]\cdot\P[E_\rho] + \E\left[\theta(0, y) \mid E_\rho^c\right]\cdot(1-\P[E_\rho])\nonumber\\
        &= \E\left[\theta^\rho(0, y)\right] + O\left(\rho|y|(1-\P[E_\rho]) + {(\rho|y|)}^{d+1}\exp(-C\rho)\right)\label{almost-done-bound},
    \end{align}
    where the last part of the last line comes from the controllability bound in Theorem~\ref{cont-estimate}. Since $\rho > C\log|y|$, we can ensure that the error term in~(\ref{almost-done-bound}) is at most $\frac{1}{2}\lambda$ and can therefore be absorbed into the constant.
\end{proof}

\subsection{Nonrandom scaling bias}
Given the controllability estimate in Theorem~\ref{cont-estimate}, we can apply Fekete's lemma~\cite{Fekete} to extract a limit, for every $v \in \R^d$, of \[ \overline{\theta}(v) := \lim_{\alpha \to \infty} \frac{\E[\theta(0, \alpha v)]}{\alpha}. \] We are now interested in bounding the nonrandom scaling bias, given by \[ \E[\theta(0, v)] - \overline{\theta}(v), \] for $v \in \R^d$. The controllability estimate implies that $\E[\theta(0, y)]-\E[\theta(0, z)] \leq C$ if $|y-z| \leq 1$, so we lose nothing by assuming that $v \in \Z^d$.

First, we find some initial sublinear bound for the nonrandom scaling bias. Then, we follow an argument of Alexander~\cite{Alexander97} to improve the bound so it matches the estimate (up to a log-factor) for the random fluctuations from Proposition~\ref{random-fluc}. We need the initial bound to estimate the region where we apply the controllability bound in the later argument.

Since $\theta(x, y)$ is subadditive, we know that $\E[\theta(0, v)] \geq \E[\overline{\theta}(v)]$. We're looking for a bound of the form \[\E[\theta(0, v)] \leq \E[\overline{\theta}(v)] + \text{error}(|v|). \] Our strategy will be to take a controlled path from $0$ to $v$, chop it into pieces, and rearrange the pieces into a new path which stays close to the line connecting $0$ to $v$. If the random fluctuations are small, the rearranged path won't take much more time than the original path. In this way, we show that \[ \E[\theta(0, nv)] + \E[\theta(0, mv)] \approx \E[\theta(0, (n+m)v)] \] for $n, m \in \N$, which will give us the desired bound.

First, we need a lemma about the rearrangement which is purely combinatorial; we include a proof based on a nearly identical result from Matoušek~\cite{Matousek}, given originally in Grinberg~\cite{Grinberg}.
\begin{lem}\label{rearrangement-lemma}
    Let $v_1, \dots v_n \in B_1$ be vectors lying in the unit ball $B_1 \subseteq \R^d$ with $\sum_{i=1}^n v_i = nx$. Then there is a permutation $\sigma \colon \{1, \dots, n\} \to \{1, \dots, n\}$ such that $\left|\sum_{i=1}^k v_{\sigma(i)} - kx\right| \leq 2d$ for every $1 \leq k \leq n$.
\end{lem}
\begin{proof}
    We say that a set $\{v_1, \dots, v_n\}$ of vectors is \textit{good} if there are coefficients $\alpha_1, \dots, \alpha_n$ such that $\alpha_i \in [0, 1]$, with \[ \sum_{i=1}^n \alpha_i v_i = (n-d)x \] and \[ \sum_{i=1}^n \alpha_i = n-d. \] Note that our original set of vectors is good, and that if a set of vectors is good then \[ \left|\sum_{i=1}^n v_i - nx\right| \leq 2d. \] To prove the lemma, it is enough to show that if $n > d$ then there is $i$ such that $\{v_1, \dots, v_n\} \setminus \{v_i\}$ is also good (then we build $\sigma$ by putting $\sigma(n) = i$ and proceeding recursively). Consider the following system of equations in $n$ unknowns $x_1, \dots, x_n$: \[ \sum_{i=1}^n x_i v_i = (n-d-1)x, \qquad \sum_{i=1}^n x_i = n-d-1. \] Note that $x_i = \alpha_i(n-d-1)/(n-d)$ shows that there is a solution with $x_i \in [0, 1]$ for all $i$. Since there are $d+1$ equations, we can modify $x_1, \dots, x_n$ so that all but $n-d-1$ of them are either $0$ or $1$. If none of the $x_i$ were $0$, then $n-d-1$ of them would $1$, so we would have $\sum_{i=1}^n x_i > n-d-1$, a contradiction. To conclude, we take $i$ such that $x_i = 0$ and now the remaining $x_i$ witness the fact that $\{v_1, \dots, v_n\} \setminus \{v_i\}$ is good.
\end{proof}

Next, we apply Lemma~\ref{rearrangement-lemma} together with the fluctuation bound in Proposition~\ref{random-fluc} to show that $\E[\theta(0, \cdot)]$ is approximately additive.
\begin{lem}\label{rearrangement-applied}
    There is a constant $C = C(d, L) > 0$ such that for any $n, m \in \N$ and $v \in \R^d$ with $|v| \geq 2$, we have \[ \E[\theta(0, (n+m)v)] \leq \E[\theta(0, nv)] + \E[\theta(0, mv)] \leq \E[\theta(0, (n+m)v)] + C|v|^{2/3}\log^3|v|. \]
\end{lem}
\begin{proof}
    Let $E$ denote the event that $|\theta(x, y) - \E[\theta(x, y)]| \leq C|v|^{1/3}\log^3|v|$ for all $x, y \in B(0, C(n+m)v)$ with $|x-y| \leq C|v|^{2/3}$. If we choose $C$ large enough, Proposition~\ref{random-fluc} implies that $E$ has probability at least $\frac12$. Working in $E$, let $\gamma \colon [0, \theta(0, (n+m)v)] \to \R^d$ be a controlled path from $0$ to $(n+m)v$. By chopping $\gamma$ into pieces of length at most $|v|^{2/3}$, we find a sequence of points $0 = x_0, x_1, \dots, x_k = (n+m)v$ with $k \leq C|v|^{1/3}$ which lie in $\Z^d$ with \[ \sum_{i=1}^k \E[\theta(x_{i-1}, x_i)] \leq \E[\theta(0, (n+m)v)] + C|v|^{2/3}\log^3|v|. \] Using $\Z^d$-translation invariance and rearranging using Lemma~\ref{rearrangement-lemma}, we assume that $|x_i-\frac{i}{k}(n+m)v| \leq C|v|^{2/3}$. So, letting $p := \lfloor nk/(n+m) \rfloor$, we have
    \begin{align*}
        \E[\theta(0, nv)] + \E[\theta(0, mv)] &\leq \sum_{i=1}^p \theta(x_{i-1}, x_i) + C|v|^{2/3} + \sum_{i=p+1}^k \theta(x_{i-1}, x_i)\\
        &\leq \E[\theta(0, (n+m)v)] + C|v|^{2/3}\log^3|v|,
    \end{align*}
    where the first inequality comes from considering the path that goes from $0$ to $x_p$, takes a detour to $nv$, then continues from $x_p$ to $x_k$.
\end{proof}

We finish by using Lemma~\ref{rearrangement-applied} and subadditivity to bound the nonrandom scaling bias.
\begin{prop}\label{initial-sublinear-bound}
    There is a constant $C = C(d, L) > 0$ such that $|\E[\theta(0, v)] - \overline{\theta}(v)| \leq C|v|^{2/3}\log^3|v|$ for all $|v| \geq 2$.
\end{prop}
\begin{proof}
    Using Lemma~\ref{rearrangement-applied}, we note that $\E[\theta(0, v)] + C|v|^{2/3}\log^3|v|$ is superadditive in $v$. Since $C|v|^{2/3}\log^3|v|$ is strictly sublinear in $|v|$, we conclude that \[ \E[\theta(0, v)] + C|v|^{2/3}\log^3|v| \geq \lim_{n \to \infty} n^{-1}\left(\E[\theta(0, nv)] + C|nv|^{2/3}\log^3|nv|\right) = \overline{\theta}(v), \] and the conclusion follows.
\end{proof}

Now that we have some initial sublinear bound on the scaling bias, we will follow Alexander's~\cite{Alexander97} argument to show that $\E[\theta(0, v)]$ satisfies the convex hull approximation property (defined below) with exponent $\frac12$. The argument is a bit different in our setting because (i) the controllability estimate only holds with high probability (not almost surely) and (ii) we don't have access to the van~den~Berg\textendash{}Kesten inequality, so we replace it with the finite range of dependence assumption together with the controllability estimate. The speed limit $L+1$ on controlled paths provides the necessary locality to apply the finite range of depencence assumption.

Let $f \colon \Z^d \to \R$ be a nonnegative subadditive function with sublinear growth; that is, there is some constant $r > 0$ such that $f(x) \leq r|x|$ (in practice we will set $f(x) := \E[\theta(0, x)]$). Define $\overline{f} \colon \R^d \to \R$ by \[ \overline{f}(x) := \lim_{n \to \infty} n^{-1}f([nx]), \] where $[\cdot]$ denotes coordinate-wise rounding to integers. Then $\overline{f}$ is positively homogeneous and convex, so we define $\overline{f}_x$ to be a supporting affine functional to $\overline{f}$ at $x$, chosen consistently so that $\overline{f}_x = \overline{f}_{\alpha x}$ for all $\alpha > 0$. We can think of $\overline{f}_x(v)$ as representing the progress that an increment $v$ makes in the direction of $x$. Since $\overline{f}$ is positively homogenous, we have $\overline{f}_x(0) = 0$ and $\overline{f}_x(x) = f(x)$.

In the following, $\varphi \colon (1, \infty) \to \R$ is a positive nondecreasing function (in practice we'll set $\varphi(x) := \log^3 x$).

\begin{defn}
    We say that $f$ satisfies the general approximation property with exponent $\nu \geq 0$ and correction factor $\varphi$ if there are constants $C, M > 0$ such that if $|x| \geq M$ then \[ \overline{f}(x) \leq f(x) \leq \overline{f}(x) + C|x|^\nu\varphi(|x|). \]
\end{defn}

Next, we define the set of ``good'' increments toward $x$ by
\begin{equation}\label{good-definition}
    G_x(\nu, \varphi, C, K) := \left\{v \in \Q^d \mid |v| \leq K|x|, \, \overline{f}_x(v) \leq \overline{f}_x(x), \, f(v) \leq \overline{f}_x(v) + C|x|^\nu\varphi(|x|) \right\},
\end{equation}
where $C, K \geq 0$. The conditions in~(\ref{good-definition}) say that good increments are not too much larger than $|x|$, don't overshoot in the direction of $x$, and aren't too wasteful along the way (we can think of $f(v)-\overline{f}_x(v)$ as the inefficiency in a step $v$ towards $x$). If we could write every $x$ as the sum of a bounded number of good increments, this would imply the general approximation property. Unfortunately, this is not so easy to accomplish. Instead, we approximate $nx$ by $O(n)$-many good increments for some large $n$; this is the content of the next definitions. In the proofs, we will need to apply the controllability estimate along a path from $0$ to $nx$, so we need some upper bound on $n$. The sublinear bound in Proposition~\ref{initial-sublinear-bound} allows us to bound $n$ by some power of $|x|$, so we can apply controllability.

\begin{defn}
    A sequence $v_0, v_1, \dots, v_m$ is a $G_x(\nu, \varphi, C, K)$-skeleton if $v_i-v_{i-1} \in G_x(\nu, \varphi, C, K)$ for all $1 \leq i \leq m$.
\end{defn}
\begin{defn}
    A function $f$ satisfies the skeleton approximation property with exponent $\nu \geq 0$ and correction factor $\varphi$ if there are constants $M, C, K > 0$ and $a > 1$ such that, for every $x \in \Q^d$ with $|x| \geq M$, there exists $n \geq 1$ and a $G_x(\nu, \varphi, C, K)$-skeleton $0 = v_0, v_1, \dots, v_m = nx$ of length $m \leq an$.
\end{defn}

As we will see, the skeleton approximation property is more convenient to verify. The following theorem of Alexander~\cite{Alexander97}, whose proof is in the appendix, provides a link between the skeleton and general approximation properties.
\begin{thm}[Alexander~\cite{Alexander97}]\label{skeleton-implies-gap}
    Suppose that $f \colon \Z^d \to \R$ is a nonnegative subadditive function and there is a constant $r \geq 1$ such that $f(x) \leq r|x|$ for all $x \in \Z^d$. If $f$ satisfies the skeleton approximation property with exponent $\nu > 0$ and correction factor $\varphi$, then $f$ satisfies the general approximation property with the same exponent and correction factor.
\end{thm}

It remains to verify that $f(x) := \E[\theta(0, x)]$ satisfies the skeleton approximation property.
\begin{prop}\label{skeleton-exists}
    The function $f(x) := \E[\theta(0, x)]$ satisfies the skeleton approximation property with exponent $\nu = \frac12$ and correction factor $\varphi(x) = \log^3|x|$.
\end{prop}
\begin{proof}
    Given a skeleton ${\{v_i\}}_{0 \leq i \leq m}$, we define its error to be \[ \err\left({\{v_i\}}_{0 \leq i \leq m}\right) := \sum_{i=0}^{m-1} \max\left(0, \E[\theta(v_i, v_{i+1})] - \theta(v_i, v_{i+1})\right), \] measuring how much faster a controlled path can traverse the skeleton than expected.

    We say that a skeleton ${\{v_i\}}_{0 \leq i \leq m}$ is $\eta$-reasonable if, for every $s > 0$ and $1 \leq i \leq j \leq m$, if $s \leq |v_i-v_{i-1}|, |v_j-v_{j-1}| \leq 2s$ and there are at least $\eta-1$ other indices $i < k < j$ such that if $s \leq |v_k-v_{k-1}| \leq 2s$, then \[ (L+1)(\E|v_i-v_{i-1}| + \E|v_j-v_{j-1}|) + 1 \leq |v_i-v_j|. \] This definition means that legs of a reasonable skeleton contribute terms to the error which are independent random variables, as long as they're of the same scale and far enough apart.

    \textit{Step 1.} We show that, with probability which tends to $1$ as $|x| \to \infty$, there is $\eta > 0$ such that every $\eta$-reasonable $G_x$-skeleton has small error. Let ${\{v_i\}}_{0 \leq i \leq m}$ be an $\eta$-reasonable $G_x$-skeleton. Then partition the indices $1, 2, \dots, m$ into $O(\log |x|)$ buckets numbered starting at $0$, where the $i$th bucket is given by \[ B_i = \{j \mid 2^i \leq |v_j-v_{j-1}| < 2^{i+1}\}. \] Within a single bucket, indices which are at least $\eta$ apart (in sorted order) contribute independent terms to the error, so splitting the bucket into $\eta$ sub-buckets containing independent terms, exponentiating the random fluctuation bound from Proposition~\ref{random-fluc}, and using Markov's inequality yields \[ \P\left[\sum_{j \in B_i} \max\left(0, \E[\theta(v_j, v_{j+1})] - \theta(v_j, v_{j+1})\right) > \zeta m \eta |x|^{1/2} \log^2 |x|\right] \leq \eta{\left(C\exp(-C^{-1}\zeta\log^2|x|)\right)}^m \] for some constant $C > 0$ and any $\zeta \geq 1$. Summing over the buckets, we get
    \begin{align}
        \P\left[\sum_{j=1}^m \max\left(0, \E[\theta(v_j, v_{j+1})] - \theta(v_j, v_{j+1})\right) > \zeta m \eta |x|^{1/2} \log^3 |x|\right] &\leq \eta\log|x|{\left(C\exp(-C^{-1}\zeta\log^2|x|)\right)}^m\nonumber\\
        &\leq C\exp(-\zeta C^{-1}m\log^2|x|),\label{bucket-bound}
    \end{align}
    where the last inequality holds as long as $|x|$ is sufficiently large.

    Estimate~(\ref{bucket-bound}) holds for a particular skeleton ${\{v_i\}}_{0 \leq i \leq m}$. Since there are at most $C|x|^{dm}$ many $G_x$-skeletons of length $m$, we can set $\zeta$ large enough to sum over all $G_x$ skeletons and get \[ \P\left[\text{every $\eta$-reasonable $G_x$-skeleton of length $m$ has error at most $Cm|x|^{1/2}\log^3|x|$}\right] \geq 1-C\exp(-C^{-1}\log^2 |x|). \] Here we make a note to choose $M$ large enough so that $C\exp(-C^{-1}\log^2 M) < \frac13$.

    \textit{Step 2.} Choose $n \in \N$ such that $\frac{1}{n}f(nx) - \overline{f}(x) < 1$. In view of Proposition~\ref{initial-sublinear-bound}, we can ensure that $n \leq C|x|^3$. Note that there always exists a $G_x$-skeleton from $0$ to $nx$, since sufficiently short increments are in $G_x$. The challenge is to find a $G_x$-skeleton without too many vertices.

    We show that, with probability which tends to $1$ as $|x| \to \infty$, there exists a reasonable $G_x$-skeleton of length $m \leq n$ from $0$ to $nx$. Proposition~\ref{random-fluc} applied to $\theta(0, nx)$ shows that \[ \P\left[\theta(0, nx) - \E[\theta(0, nx)] > n\right] < \frac13, \] as long as $M$ (and therefore $|x|$) is large enough. Similarly, we can choose $M$ large enough to ensure that \[ \P\left[\mathcal{E}\left(2(L+1)\E[\theta(0, nx)]\right) > |x|^{1/4}\right] < \frac13. \] By the union bound, there is some $\omega \in \Omega$ such that \[ \theta_\omega(0, nx) \leq \E[\theta(0, nx)] + n \] and
    \begin{equation}\label{control-bound-subadditive}
        \mathcal{E}_\omega(2(L+1)\E[\theta(0, nx)]) \leq |x|^{1/4}
    \end{equation}
    and so that for every $m \geq 1$, every $G_x$-skeleton of length $m$ has error no more than $Cm|x|^{1/2}\log^3|x|$.

    We now build a $G_x$-skeleton greedily from $\gamma \colon [0, \theta_\omega(0, nx)] \to \R^d$, the shortest controlled path from $0$ to $nx$ with respect to this particular $\omega$. Let $v_0 := 0$ and $t_0 := 0$. Given $v_i, t_i$ for $i \geq 0$, define \[ t_{i+1} := \min \left(\theta_\omega(0, nx), \sup \{ t \geq t_i \mid (\gamma(t)-v_i + B_1) \cap G_x \neq \emptyset \}\right). \] Then, define $v_{i+1}$ to be any point in $(\gamma(t_{i+1}) + \overline{B_1}) \cap (G_x + v_i)$. Continue until $v_{i+1} = nx$, and set $m := i+1$ in this case.

    We claim that the skeleton $v_0, \dots, v_m$ is $\eta$-reasonable for some constant $\eta = \eta(d, L) > 0$. Indeed, suppose not. Then there would be some $s > 0$ and indices $i_0 < i_1 < \cdots < i_\eta$ such that $s \leq |v_{i_{j+1}} - v_{i_j}| \leq 2s$ for all $0 \leq j < \eta$ and $|v_{i_\eta}-v_{i_0}| \leq Cs + 1$, where $C(d, L) > 0$ is a constant; for example we could take $C = 4(L+1)\sup_{|v|=1}\E[\theta(0, v)]$, which depends only on $d$ and $L$ by Theorem~\ref{cont-estimate}. Since we chose the skeleton greedily, we can be sure that
    \begin{equation}\label{step-big}
        s \geq \frac12 |v_{i+1}-v_i| \geq C^{-1}{|x|}^{1/2},
    \end{equation}
    where we may have to enlarge the constant $C > 0$. Define $\tilde{\gamma}$ to be the path which follows $\gamma$ from $0$ to $\gamma(t_{i_0})$, then follows a shortest controlled path from $\gamma(t_{i_0})$ to $\gamma(t_{i_\eta})$, then follows $\gamma$ the rest of the way to $nx$. The controllability bound~(\ref{control-bound-subadditive}) implies that \[ \len(\tilde{\gamma}) \leq \len(\gamma) - \eta {(L+1)}^{-1} s + Cs + 1 + |x|^{1/4}. \] From~(\ref{step-big}) we see that $\len(\tilde{\gamma})$ is strictly smaller than $\len(\gamma)$ if $\eta \geq 4C(L+1)$ (here we make a note to choose $C,M \geq 1$). This contradicts the fact that $\gamma$ is a shortest path, so we conclude that the skeleton is $\eta$-reasonable.

    Now we use the error bound for reasonable skeletons to see that \[ \err\left({\{v_i\}}_{1 \leq i \leq m}\right) \leq C m|x|^{1/2}\log^3|x|, \] and so using the fact that $\gamma$ is a shortest path we see that
    \begin{align}
        \sum_{i=0}^{m-1} \E[\theta_\omega(v_i, v_{i+1})] &\leq \theta_\omega(0, nx) + m(4 + |x|^{1/4}) + Cm|x|^{1/2}\log^3|x| \nonumber\\
        &\leq \E[\theta_\omega(0, nx)] + n + m(4 + |x|^{1/4}) + Cm|x|^{1/2}\log^3|x|. \label{skeleton-short-bound}
    \end{align}

    On the other hand, we can split the indices $0, 1, 2, \dots, m-1$ into two groups: define \[ L := \{0 \leq i < m \mid |v_{i+1}-v_i| + \sqrt{d} \geq K|x| \text{ or } \overline{f}_x(v_{i+1}-v_i) + r\sqrt{d} \geq \overline{f}(x)\}, \] and let $S := \{0, \dots, m-1\} \setminus L$. We think of $S$ as the indices of short increments and $L$ as the indices of long increments. Note that every index $i \in S$ of a short increment satisfies \[ \E[\theta_\omega(v_i, v_{i+1})] - \overline{f}_x(v_{i+1}-v_i) \geq C|x|^{1/2}\log^3|x| - O(1). \] So summing over $i \in S$ and choosing $C > 0$ large enough relative to the bound~(\ref{skeleton-short-bound}) shows that we can ensure $|S| \leq \frac{m}{4}$. On the other hand, the linear growth of $\overline{f}$ shows that there are at most $C|x|^{-1}\theta_\omega(0, nx) \leq Cn$ long increments, so $m = O(n)$ as desired.
\end{proof}

\subsection{A shape theorem}
To conclude the section, we put together our bounds on random and nonrandom error to get a quantitative shape theorem for the metric problem. Define the large-scale reachable set at time $t$ by \[ \mathcal{S}_t := \{x \in \R^d \mid \overline{\theta}(x) \leq t\}. \]
\begin{thm}\label{shape-theorem}
    Let $V \colon \R^d \to \R^d$ be a random divergence-free vector field with unit range of dependence and $\|V\|_{C^{1,1}} \leq L$ almost surely. Then there are constants $C(d, L) > 1 > c(d, L) > 0$ such that, for all $t \geq 0$, \[ \P\left[\dist_H(\mathcal{R}_t(0), \mathcal{S}_t) > Ct^{1/2}\log^3 t + \lambda\right] \leq C\exp\left(\frac{-C^{-1}\lambda^{1/2}}{t^{1/4}}\right), \] where $\dist_H$ denotes the Hausdorff distance. Furthermore, there is a random variable $T_0$, with \[ \E[\exp(c\log^{3/2} T_0)] < \infty, \] such that \[ \sup_{(t,x) \in [0, T] \times B_T} \frac{\dist_H(\mathcal{R}_t(x), x + \mathcal{S}_t)}{T^{1/2}\log^3 T} \leq C \] for all $T \geq T_0$.
\end{thm}
\begin{proof}
    For the first claim, let $t \geq 0$. Apply Theorem~\ref{cont-estimate} to $B_{(L+1)t + 1}$ and Proposition~\ref{random-fluc} to every $x \in \Z^d \cap \overline{B_{(L+1)t}}$ and use the union bound to see that as long as $\lambda \geq Ct^{1/2}\log^2 t$ we have
    \begin{equation}\label{random-fluc-applied}
        \P\left[\forall x \in \overline{B_{(L+1)t}} :\: \left|\theta(0, x) - \E[\theta(0, x)]\right| > \lambda \right] \leq C\exp\left(\frac{-C^{-1}\lambda^{1/2}}{t^{1/4}}\right),
    \end{equation}
    where we absorbed polynomials into the exponential by enlarging the constant $C$. Note also that Theorem~\ref{cont-estimate} implies that if $0 \leq r \leq s$, then
    \begin{equation}\label{cont-estimate-applied}
        \P\left[\mathcal{R}_r(0) \nsubseteq \mathcal{R}_s(0) + B_{\lambda}\right] \leq C\exp\left(\frac{-C^{-1}\lambda^{1/2}}{t^{1/4}}\right).
    \end{equation}
    This bounds the random error. On the other hand, Proposition~\ref{skeleton-exists} and Theorem~\ref{skeleton-implies-gap} combine to show that
    \begin{equation}\label{nonrandom-bias-applied}
        0 \leq \E[\theta(0, x)] - \overline{\theta}(x) \leq C|x|^{1/2}\log^3|x|,
    \end{equation}
    which bounds the nonrandom error. The estimates~(\ref{random-fluc-applied}) and~(\ref{nonrandom-bias-applied}) say that, with high probability, the first passage time $\theta(0, x)$ from $0$ to any point $x$ is close to the large-scale average $\overline{\theta}(x)$. Furthermore, the estimate~(\ref{cont-estimate-applied}) says that once a controlled path reaches $x$, the reachable set stays close to $x$ for all later times (the controllability estimate guarantees the existance of controlled paths in the form of short loops). Unwrapping the definition of Hausdorff distance yields the first claim.

    For the second claim, apply the first claim to every $(t, x) \in (\Z \cap [0, T]) \times (\Z^d \cap B_T)$ and the union bound to conclude that \[ \P\left[\sup_{t \in \Z \cap [0, T]} \sup_{x \in \Z^d \cap B_T} \dist_H(\mathcal{R}_t(x), x + \mathcal{S}_t) > CT^{1/2}\log^3 T + \lambda\right] \leq CT^{d+1}\exp\left(\frac{-C^{-1}\lambda^{1/2}}{T^{1/4}}\right). \] Next, apply the controllability estimate in $B_T$ to see that the same holds for all $(t, x) \in [0, T] \times B_T$, by enlarging the constant $C$. Plugging in $\lambda = CT^{1/2}\log^3 T$ shows that \[ \P\left[\sup_{(t, x) \in [0, T] \times B_T} \frac{\dist_H(\mathcal{R}_t(x), x+\mathcal{S}_t)}{T^{1/2}\log^3 T} > C\right] \leq C\exp(-C^{-1}\log^{3/2} T), \] and the conclusion follows.
\end{proof}

\section{Proofs of the main results}\label{applications-section}
\subsection{Homogenization of solutions}
Now that we have a quantitative shape theorem for the metric problem, we use the control formulation to extract a rate of convergence for solutions of the G equation to the large-scale limit. Let $u_0 \in C^{0,1}(\R^d)$ and let $\overline{u} \colon \R_t \times \R^d_x \to \R$ be a solution to the effective problem~(\ref{macro-g-eqn}), where the effective Hamiltonian is given by
    \begin{equation}\label{effective-H-formula}
        \overline{H}(p) := \sup_{v \in \R^d} p \cdot \frac{v}{\overline{\theta}(v)}.
    \end{equation}
    Note that $\overline{\theta}$ is positively homogeneous, so the supremum can be restricted to the unit sphere and is in fact a maximum. We have the representation formula
    \begin{equation}
        \overline{u}(t, x) = \sup_{x + \mathcal{S}_t} u_0. \label{u-rep-form}
    \end{equation}
    On the other hand, let $u^\varepsilon$ be a solution to the G equation~(\ref{eps-g-eqn}). We have the representation formula
    \begin{equation}
        u^\varepsilon(t, x) = \sup_{\varepsilon \mathcal{R}_{\varepsilon^{-1}t}(\varepsilon^{-1}x)} u_0. \label{u-epsilon-rep-form}
    \end{equation}

\begin{proof}[Proof of Theorem~\ref{main-homog-theorem}]
    Using the representation formulas~(\ref{u-rep-form}) and~(\ref{u-epsilon-rep-form}), we see that for every $0 \leq t \leq T$ and $x \in B_T$ we have \[ |u^\varepsilon(t, x) - \overline{u}(t, x)| = \left|\sup_{\varepsilon \mathcal{R}_{\varepsilon^{-1}t}(\varepsilon^{-1}x)} u_0 - \sup_{x + \mathcal{S}_t} u_0\right| \leq \|u_0\|_{C^{0,1}}\dist_H(\varepsilon\mathcal{R}_{\varepsilon^{-1}t}(\varepsilon^{-1}x), x + \mathcal{S}_t). \] Rescaling by $\varepsilon^{-1}$ and applying Theorem~\ref{shape-theorem} (using the fact that $\varepsilon^{-1}\mathcal{S}_t = \mathcal{S}_{\varepsilon^{-1}t}$) yields \[ \sup_{(t, x) \in [0, T] \times B_T} \dist_H(\varepsilon\mathcal{R}_{\varepsilon^{-1}t}(\varepsilon^{-1}x), x + \mathcal{S}_t) \leq C{(T\varepsilon)}^{1/2}\log^3 (\varepsilon^{-1}T) \] for all $T \geq \varepsilon T_0$. Combining these gives \[ |u^{\varepsilon}(t, x) - \overline{u}(t, x)| \leq C\|u_0\|_{C^{0,1}}{(t\varepsilon)}^{1/2}\log^3(\varepsilon^{-1}t) \] for all $t \geq \varepsilon T_0$ as desired.
\end{proof}

\subsection{Continuity of the effective Hamiltonian}
Now we apply Theorem~\ref{shape-theorem} to show that the effective Hamiltonian $\overline{H}$ depends continuously on the law of the environment.

\begin{proof}[Proof of Theorem~\ref{continuous-dependence-theorem}]
    Throughout this proof we will use the superscript $n$ to denote the corresponding object (Hamiltonian, first passage time, reachable set, etc.) for the environment $\P^n$. Since the effective Hamiltonians are positively homogeneous, it suffices to show uniform convergence on $B_1$. From the formula~(\ref{effective-H-formula}) for the effective Hamiltonian, we deduce that \[ \overline{H}(p) = \sup_{v \in \mathcal{S}_1} v \cdot p \] and \[ \overline{H}^n(p) = \sup_{v \in \mathcal{S}_1^n} v \cdot p. \] So, it would suffice to show that $\dist_H(\mathcal{S}_1^n, \mathcal{S}_1) \to 0$ as $n \to \infty$, where $\dist_H$ denotes the Hausdorff distance. From Theorem~\ref{shape-theorem}, we have \[ \E^n[\dist_H(t^{-1}\mathcal{R}_t(0), \mathcal{S}_1^n)] \leq Ct^{-1/2}\log^3 t \] for all $n \in \N$ and $t \geq 0$. On the other hand, since $\dist_H(t^{-1}\mathcal{R}_t(0), \mathcal{S}_1)$ is a continuous function on $L^\infty(B_{Ct}; \R^d)$, we use weak convergence to see that \[ \lim_{n \to \infty} \E^n[\dist_H(t^{-1}\mathcal{R}_t(0), \mathcal{S}_1)] = \E[\dist_H(t^{-1}\mathcal{R}_t(0), \mathcal{S}_1)] \leq Ct^{-1/2}\log^3 t. \] Combining these with the triangle inequality, we have \[ \dist_H(\mathcal{S}_1^n, \mathcal{S}_1) \leq Ct^{-1/2}\log^3 t. \] Send $t \to 0^+$ to conclude.
\end{proof}

\printbibliography{}

\begin{appendices}
    \section{Proof of Theorem~\ref{skeleton-implies-gap}}
    \textit{Step 1.} We show that if $x \in \Q^d$ with $|x| \geq C$, then there is $\alpha \in [c, 1]$ such that $\alpha x$ lies in the convex hull of $G_x$. Let $n \in \N$ be large enough so that \[ |n^{-1}f(nx) - \overline{f}(x)| \leq 1. \] Let $x_0, x_1, \dots, x_m$ be an $G_x$-skeleton for $nx$. Then \[ x = \frac{1}{n}\sum_{k=1}^m (x_k-x_{k-1}), \] and $n \leq m \leq Cn$, where the first part of the inequality follows from applying $\overline{f}_x$ to both sides of the equation.

    \textit{Step 2.} We show that if $x \in \Q^d$, $|x| \geq K$, $t \geq 1$, and $tx \in \Z^d$, then there is a $z \in \Z^d$ with $f(tx) - \overline{f}(tx) \leq f(z) - \overline{f}(z) + tC|x|^\nu\varphi(|x|)$. Using the previous step, write $tx = z + \sum_{k=1}^m v_k$, where $|z| \leq C|x|$, $\overline{f}(z) \leq \overline{f}_x(z) + C$, $v_k \in G_x$, and $m \leq Ct$. Indeed, for some $\alpha \in [c, 1]$ we first write \[ \alpha x = \sum_{i=1}^{d+1} p_i v_i, \] where $v_i \in G_x$ and $p_i \geq 0$, $\sum_i p_i = 1$. Note that the sum only requires $d+1$ terms by Caratheodory's theorem on convex hulls, since we are working in $\R^d$. To decompose $tx$, we write \[ tx = \sum_{i=1}^{d+1}(t\alpha^{-1}p_i - \lfloor t\alpha^{-1}p_i \rfloor)v_i + \sum_{i=1}^{d+1}\lfloor t\alpha^{-1}p_i \rfloor v_i =: z + (tx-z), \] so $z$ satisfies the required properties. By subadditivity of $f$ and linearity of $\overline{f}_x$,
    \begin{align*}
        f(tx) &\leq f(z) + \sum_{k=1}^m f(v_k)\\
        &\leq f(z) + \sum_{k=1}^m \overline{f}_x(v_k) + C|x|^\nu\varphi(|x|)\\
        &= f(z) + \overline{f}_x(tx - z) + tC|x|^\nu\varphi(|x|).
    \end{align*}
    Finally, we write $\overline{f}(tx) = \overline{f}_x(z) + \overline{f}_x(tx-z)$ and subtract from both sides of the inequality above to get \[ f(tx) - \overline{f}(tx) \leq f(z) - \overline{f}(z) + Ct, \] where we used the fact that $\overline{f}(z) \leq \overline{f_x}(z) + C$.

    \textit{Step 3.} For some large $M > 1$ (and possibly enlarged $K$), the previous step yields
    \begin{align*}
        \sup_{|x| \leq M^{k+1}K} f(x)-\overline{f}(x)  &\leq \sup_{|x| \leq M^k K} f(x)-\overline{f}(x) + CM|x|^\nu\varphi(|x|)\\
        &\leq \sup_{|x| \leq M^k K} f(x)-\overline{f}(x) + C(M^\nu-1)|x|^\nu\varphi(|x|),
    \end{align*}
    where we made $C$ larger in the second line to account for the change in constant. By induction on $k$, we conclude that \[ f(x)-\overline{f}(x) \leq C|x|^\nu\varphi(|x|). \]
\end{appendices}

\end{document}